\documentclass[12pt,leqno]{article}
\title{Singular Integers and Kummer-Stickelberger relation }
\author{Roland Qu\^eme}
\usepackage{amsmath,amsthm,amstext,amsbsy}
\newtheorem{thm}{Theorem}[section]
\newtheorem{cor}[thm]{Corollary}

\newtheorem{lem}[thm]{Lemma}
\font\mathbb=msbm10
\newcommand{\N}{\mbox{\mathbb N}}

\newcommand{\Q}{\mbox{\mathbb Q}}
\newcommand{\Z}{\mbox{\mathbb Z}}

\newcommand{\modu}{\ \mbox{mod}\ }
\newcommand{\be}{\begin{equation}}
\newcommand{\ee}{\end{equation}}
\newcommand{\bd}{\begin{displaymath}}
\newcommand{\ed}{\end{displaymath}}
\newcommand{\bn}{\begin{enumerate}}
\newcommand{\en}{\end{enumerate}}

\date{2007 feb 07}
\begin{document}
\maketitle
\tableofcontents
\abstract
\abstract Let $p$ be an odd prime. Let ${\bf F}_p$ be the finite
field of $p$ elements and ${\bf F}_p^*$ its multiplicative group.
Let $K=\Q(\zeta)$ be a $p$-cyclotomic field and $O_K$  its ring of
integers. Let $\pi$ be the prime ideal of $K$ lying over $p$. Let
$\sigma :\zeta\rightarrow\zeta^v$ be the $\Q$-isomorphism of $K$ for
a primitive root $v\modu p$. The subgroup  of exponent $p$ of the
class group of $K$ can be seen as a direct sum $C_p=\oplus_{i=1}^r
\Gamma_i$ of groups of order $p$ where the group $\Gamma_i$ are
annihilated by  polynomials $\sigma-\mu_i$ with $\mu_i\in{\bf
F}_p^*$. Let us fix  $\Gamma$ for one  of these $\Gamma_i$. From
Kummer,   there exist non-principal prime
ideals $\mathbf q$ of inertial degree $1$ with classes $Cl(\mathbf q)\in\Gamma$. We show that there exist  singular
semi-primary integers $A$ such that $A O_K= \mathbf q^p$ and $A^{\sigma-\mu}=\alpha^p$ with $\alpha\in K$.
Let $E$ defined by $E=(\frac{\alpha}{\overline{\alpha}})^p-1$.
Let $\nu$ be the
positive integer defined by $\nu=v_\pi(E)-(p-1)$ where $v_\pi(.)$ is the $\pi$-adic valuation.
The aims of this article are:
\bn
\item
To give some algebraic relations connecting the singular integer  $A$ and the $p^2$-power of the Gauss sum $g(\mathbf q)$.
\item
to give some informations on   the $\pi$-adic expansion of $g(\mathbf q)^{p^2}$ and of $E$.
\item
to derive an explicitly computable upper bound
of $\nu$ from the algebraic expression of Jacobi resolvents  used in Kummer-Stickelberger relation of the
 $p$-cyclotomic field $K$.
\en
 {\it Remark:} we have not found in the literature some formulations corresponding to
 lemma \ref{l602061} p. \pageref{l602061}, theorems \ref{t601311} p.\pageref{t601311}, \ref{l512164} p. \pageref{l512164},
 \ref{t608201} p. \pageref{t608201}, \ref{t605031} p. \pageref{t605031} and \ref{t602061} p.  \pageref{t602061}.
%
%
\section{Some definitions}\label{s601191}
In this section we give the definitions and  notations on cyclotomic
fields,  $p$-class group, singular  numbers,  primary and non-primary, used in this paper.
\begin{enumerate}
\item
Let $p$ be an odd prime. Let $\zeta_p$ be a root of the polynomial equation $X^{p-1}+X^{p-2}+\dots+X+1=0$.
Let $K_p$ be the $p$-cyclotomic field $K_p=\Q(\zeta_p)$. The ring of integers of $K_p$ is $\Z[\zeta_p]$.
Let $K_p^+$ be the maximal totally real subfield of $K_p$.
The ring of integers of $K_p^+$ is $\Z[\zeta_p+\zeta_p^{-1}]$ with group of units $\Z[\zeta_p+\zeta_p^{-1}]^*$.
Let $v$ be a primitive root $\modu p$ and $\sigma: \zeta_p\rightarrow \zeta_p^v$ be a $\Q$-isomorphism of $K_p$.
In this article, for $n\in\Z$, following the conventions adopted by Ribenboim \cite{rib},
we use the notation $v_n=v^n\modu p$ with $1\leq v_n\leq p-1$.
Let $G_p$ be the Galois group of the extension $K_p/\Q$.
Let ${\bf F}_p$ be the finite field of cardinal $p$ with no null part  ${\bf F}_p^*$.
Let $\lambda=\zeta_p-1$. The prime ideal of $K_p$ lying over $p$ is $\pi=\lambda \Z[\zeta_p]$.
\item
Let $C_p$ be the $p$-class group of $K_p$ (the subgroup of exponent
$p$ of the class group $\mathbf C$ of $K_p$). Let $r$ be the rank of
$C_p$. Let $C_p^+$ be the $p$-class group of $K_p^+$. Let $r^+$ be
the rank of $K_p^+$.  Then $C_p=C_p^+\oplus C_p^-$ where $C_p^-$ is
the relative $p$-class group.
\item
$C_p$ is the direct sum of $r$ subgroups  $\Gamma_i$  of order $p$
annihilated by $\sigma-\mu_i\in {\bf F}_p[G_p]$ with $\mu_i\in{\bf
F}_p^*$
\begin{equation}\label{e611031}
C_p=\oplus_{i=1}^r \Gamma_i.
\end{equation}
Then $\mu= v_{n}$ with a natural integer $n,\quad 1\leq
n\leq p-2$.
\item
An  integer  $A\in \Z[\zeta_p]$ is said singular if $A^{1/p}\not\in
K_p$ and if  there exists an   ideal $\mathbf a$  of $\Z[\zeta_p]$
such that $A \Z[\zeta_p]=\mathbf a^p$.
\item
An integer  $A\in \Z[\zeta_p]$ is said semi-primary if there exists a natural integer $a$ coprime with $p$ such that
$A\equiv a\modu\pi^2$.
An integer  $A\in \Z[\zeta_p]$ is said primary if there exists a natural integer $a$ coprime with $p$ such that
$A\equiv a^p\modu\pi^p$.
%
\section{Preliminary results}
We shall use in this paper a result of Kummer:  the group of ideal classes of the
cyclotomic field $K_p$ is generated by the classes of prime ideals
$\mathbf q$ of inertial degree $1$.
Let $\Gamma$ be one of the groups $\Gamma_i$ defined in  the relation (\ref{e611031}).
\begin{enumerate}
\item
\underline {If  $r^->0$ and  $\Gamma\subset  C_p^-$: }
then there exist  singular semi-primary integers $A$  with $A \Z[\zeta_p] =\mathbf q^{2p}$
where $\mathbf q$ is a {\bf non}-principal prime  ideal of $\Z[\zeta_p]$ verifying   simultaneously
\begin{equation}\label{e512101}
\begin{split}
& Cl(\mathbf q)\in \Gamma,\\
& \sigma(A)=A^\mu\times\alpha^p,\quad \mu\in {\bf F}_p^*,\quad \alpha\in K_p,\\
&\mu= v_{2m+1}, \quad m\in\N, \quad 1\leq m\leq \frac{p-3}{2},\\
&\pi^{2m+1} \ |\ A-a^p,\quad a\in\N,\quad 1\leq a\leq p-1,\\
\end{split}
\end{equation}
Moreover, this number $A$ verifies
\begin{equation}\label{e512103}
A\times\overline{A}=D^p,
\end{equation}
for some integer $D\in \Z[\zeta_p+\zeta_p^{-1}]^*$.
\begin{enumerate}
\item
If  $A$ is   non-primary  then $\pi^{2m+1} \ \|\ A-a^p$.
\item
If   $A$ is    primary then  the extension
$K(A^{1/p})/K$ is cyclic unramified of degree $p$. In that case we know that $r^+>0$.
\end{enumerate}
(see Qu\^eme, \cite{que} theorem 2.4 p. 4 for these results).
Observe that there exists such singular primary numbers $A$ iff Vandiver's conjecture is false and that the  knowledge of $\pi$-adic value
$v_\pi(A-a^p)$ is important in regards of Vandiver's conjecture.
\item
\underline {If $r^+>0$ and    $\Gamma\subset C_p^+$: }
then there exist  singular integers $A$  with $A \Z[\zeta_p] =\mathbf q^{2p}$
where $\mathbf q$ is a {\bf non}-principal prime ideal of $\Z[\zeta_p]$ verifying   simultaneously
\begin{equation}\label{e6012210}
\begin{split}
& Cl(\mathbf q)\in \Gamma,\\
& \sigma(A)=A^\mu\times\alpha^p,\quad \mu\in {\bf F}_p^*,\quad \alpha\in K_p,\\
&\mu= v_{2m}, \quad m\in\N, \quad 1\leq m\leq \frac{p-3}{2},\\
&\pi^{2m} \ |\ A-a^p,\quad a\in\Z,\quad 1\leq a\leq p-1,\\
\end{split}
\end{equation}
Moreover, this number $A$ verifies
\begin{equation}\label{e512103}
\frac{A}{\overline{A}}=D^p,
\end{equation}
for some number  $D\in K_p^+$.
\begin{enumerate}
\item
If  $A$ is  non-primary then $\pi^{2m} \ \|\ A-a^p$.
\item
If  $A$ is   primary then  the extension $K(A^{1/p})$ is cyclic unramified of degree $p$.
\end{enumerate}
(see Qu\^eme, \cite{que} theorem 2.7 p. 7 for these results).
\end{enumerate}
\end{enumerate}
%
%
\section{On Kummer-Stickelberger relation}\label{s601192}
Stickelberger relation was already known by Kummer under the form of
Jacobi resolvents for the case of cyclotomic field $K_p$, see for instance
Ribenboim \cite{rib} (2.6) p. 119, which explains the title of this section.
\begin{enumerate}
\item
Let $q\not=p$ be an odd prime.
Let $\zeta_q$ be a root of the minimal polynomial equation $X^{q-1}+X^{q-2}+\dots+X+1=0$.
Let $K_q=\Q(\zeta_q)$ be the $q$-cyclotomic field.
The ring of integers of $K_q$ is $\Z[\zeta_q]$.
Here we fix a notation for the sequel.
Let $u$ be a primitive root $\modu q$.
In this article, for $n\in\Z$ we use the notation $u_n=u^n\modu q$ with $1\leq u_n\leq q-1$.
Let $K_{pq}=\Q(\zeta_p,\zeta_q)$. Then $K_{pq}$ is the compositum $K_pK_q$.
The ring of integers of $K_{pq}$ is $\Z[\zeta_{pq}]$.
\item
Let $\mathbf q$ be a prime ideal of $\Z[\zeta_p]$ lying over the
prime $q$. Let $m=N_{K_p/\Q}(\mathbf q)= q^f$ where $f$ is the order
of $q\modu p$. If $\psi(\alpha)=a$
is the image of $\alpha\in \Z[\zeta_p]$ under the natural map $\psi:
\Z[\zeta_p]\rightarrow \Z[\zeta_p]/\mathbf q$, then for
$\psi(\alpha)=a\not\equiv 0$ define a character $\chi_{\mathbf
q}^{(p)}$ on ${\bf F}_m=\Z[\zeta_p]/\mathbf q$ by
\begin{equation}
\chi_{\mathbf q}^{(p)}(a)={\{\frac{\alpha}{\mathbf q}\}}_p^{-1}=\overline{\{\frac{\alpha}{\mathbf q}\}}_p,
\end{equation}
where $\{\frac{\alpha}{\mathbf q}\}=\zeta_p^c$ for some natural integer $c$,
is the $p^{th}$ power residue character $\modu \mathbf q$.
We define
\begin{equation}\label{e6012211}
g(\mathbf q)=\sum_{x\in{\bf F}_m}(\chi_{\mathbf q}^{(p)}(x)\times\zeta_q^{Tr_{{\bf F}_m/{\bf F}_q}(x)})\in \Z[\zeta_{pq}].
\end{equation}
It follows that $\mathbf g(\mathbf q)\in \Z[\zeta_{pq}]$.
Moreover it is known that $g(\mathbf q)^p\in \Z[\zeta_p]$, see for instance Mollin \cite{mol} prop. 5.88 (c) p. 308
and that $g(\mathbf q)^{\sigma-v}\in K_p$, see for instance Ribenboim \cite{rib} (2A) b. p. 118.
\end{enumerate}
%
The Stickelberger's relation used in this article is classically:
\begin{equation}\label{e512121}
\mathbf g(\mathbf q)^p O_K=\mathbf q^{S},
\end{equation}
with $S=\sum_{t=1}^{p-1} t\times \varpi_t^{-1}$
where  $\varpi_t\in Gal(K_p/\Q)$ is given by $\varpi_t: \zeta_p\rightarrow \zeta_p^t$.
See for instance Mollin \cite{mol} thm. 5.109 p. 315.
We remind that $g(\mathbf q)\in\Z[\zeta_p,\zeta_q]$, $g(\mathbf q)\overline{g(\mathbf q)}=q^f$ and $g(\mathbf q)^p\in \Z[\zeta_p]$,
see for instance Mollin \cite{mol} prop. 5.88 p. 308.
In following lemma we give another straightforward form of $S$ more adequate for computations in $\Z[G_p]$.
%
\begin{lem}\label{l12161}
Let $P(\sigma)=\sum_{i=0}^{p-2} \sigma^i\times v_{-i}\in\Z[G_p]$.
Then $S=P(\sigma)$.
\begin{proof}
Let us consider one term $\varpi_t^{-1} \times t$.
Then $v_{-1}$ is a primitive root $\modu p$, hence  there exists one and one $i$ such that
$t=v_{-i}$. Then $\varpi_{v_{-i}}:\zeta_p\rightarrow \zeta_p^{v_{-i}}$ and so $\varpi_{v_{-i}}^{-1}:\zeta_p\rightarrow\zeta_p^{v_i}$
and so $\varpi_{v_{-i}}^{-1}=\sigma^i$, which achieves the proof.
\end{proof}
\end{lem}
%
\subsection{Some properties  of $\mathbf g(\mathbf q)$.}
In this subsection we give some elementary properties of  $g(\mathbf q)$ and, when $\mathbf q$ is of degree $1$,
we derive a Kummer form of the   Stickelberger relation where $g(\mathbf q)$ is a Jacobi resolvent.
\begin{lem}\label{l512151}
If $q\not\equiv 1\modu p$ then $g(\mathbf q)\in \Z[\zeta_p]$.
\begin{proof}$ $
\begin{enumerate}
\item
Let $u$ be a primitive root $\modu q$. Let $\tau :\zeta_q\rightarrow
\zeta_q^u$ be a $\Q$-isomorphism generating $Gal(K_q/\Q)$. The
$\Q$-isomorphism $\tau$ is extended to a $K_p$-isomorphism of
$K_{pq}$ by $\tau:\zeta_q\rightarrow \zeta_q^u,\quad
\zeta_p\rightarrow \zeta_p$. Then  $g(\mathbf q)^p\in\Z[\zeta_p]$ and so
\begin{displaymath}
\tau(g(\mathbf q))^p=g(\mathbf q)^p,
\end{displaymath}
and it follows that there exists a natural integer $\rho$ with $\rho<p$ such that
\begin{displaymath}
\tau(g(\mathbf q))= \zeta_p^\rho\times  g(\mathbf q).
\end{displaymath}
Then $N_{K_{pq}/K_p}(\tau(g(\mathbf q)))=\zeta_p^{(q-1)\rho}\times N_{K_{pq}/K_p}(g(\mathbf q))$ and so  $\zeta_p^{\rho(q-1)}=1$.
\item
If $q\not\equiv 1\modu p$, it implies that $\zeta_p^\rho=1$ and so that $\tau(g(\mathbf q))=g(\mathbf q)$
and thus  $g(\mathbf q)\in \Z[\zeta_p]$.
\end{enumerate}
\end{proof}
\end{lem}
%
Let us note in the sequel $g(\mathbf q)=\sum_{i=0}^{q-2} g_i\times \zeta_q^i$ with $g_i\in \Z[\zeta_p]$.
\begin{lem}\label{l512152}
If $q\equiv 1\modu p$ then $g_0=0$.
\begin{proof}
Suppose that $g_0\not=0$ and search for a contradiction:
we start of
\begin{displaymath}
\tau(g(\mathbf q))= \zeta_p^\rho\times  g(\mathbf q).
\end{displaymath}
We have $g(\mathbf q)=\sum_{i=0}^{q-2} g_i\times \zeta_q^i$  and so
$\tau(g(\mathbf q))=\sum_{i=0}^{q-2}  g_i\times \zeta_q^{i u}$,
therefore
\begin{displaymath}
\sum_{i=0}^{q-2} (\zeta_p^\rho \times g_i)\times \zeta_q^i=\sum_{i=0}^{q-2}  g_i\times \zeta_q^{i u},
\end{displaymath}
thus $g_0=\zeta_p^\rho \times g_0$ and so $\zeta_p^\rho=1$ which
implies that $\tau(g(\mathbf q))=g(\mathbf q)$ and so $g(\mathbf q)\in \Z[\zeta_p]$.
Then in $\Z[\zeta_p]$, from  Stickelberger relation, $g(\mathbf q)^p \Z[\zeta_p] =\mathbf q^S$ where $S=\sum_{t=1}^{p-1}t\times\varpi_t^{-1}$.
Therefore $\varpi_1^{-1}(\mathbf q) \ \|\ \mathbf q^S$ because $q$ splits totally in $K_p/\Q$
 and $\varpi_t^{-1}(\mathbf q)\not=\varpi_{t^\prime}^{-1}(\mathbf q)$ for $t\not=t^\prime$.
 This case is not possible because the first member $g(\mathbf q)^p$ is a $p$-power.
\end{proof}
\end{lem}
%
Here we give an elementary computation of $g(\mathbf q)$ not involving directly the Gauss sums.
\begin{lem}\label{l512152}
If $q\equiv 1\modu p$ then
\begin{equation}\label{e512151}
\begin{split}
&g(\mathbf q)=\zeta_q +\zeta_p^\rho\zeta_q^{u_{-1}}+\zeta_p^{2\rho}\zeta_q^{u_{-2}}+\dots +\zeta_p^{(q-2)\rho}\zeta_q^{u_{-(q-2)}},\\
& g(\mathbf q)^p \Z[\zeta_p] =\mathbf q^S,\\
\end{split}
\end{equation}
for some natural number $\rho,\quad 1<\rho\leq p-1$.
\begin{proof} $ $
\begin{enumerate}
\item
We start of $\tau(g(\mathbf q))=\zeta_p^\rho\times g(\mathbf q)$ and so
\begin{equation}\label{e512152}
\sum_{i=1} ^{q-2}g_i \zeta_q^{ui}=\zeta_p^\rho\times\sum_{i=1}^{q-2} g_i \zeta_q^i,
\end{equation}
which implies that $g_i=g_1\zeta_p^\rho $ for  $u\times i\equiv 1\modu q$ and so $g_{u_{-1}}=g_1\zeta_p^\rho$ (where $u_{-1}$
is to be understood by
$u^{-1}\modu q$,  so $1\leq u_{-1}\leq q-1)$.
\item
Then
$\tau^2(g(\mathbf q))=\tau(\zeta_p^{\rho} g(\mathbf q))=\zeta_p^{2\rho} g(\mathbf q)$.
Then
\begin{displaymath}
\sum_{i=1} ^{q-2}g_i \zeta_q^{u^2i}=\zeta_p^{2\rho}\times (\sum_{i=1}^{q-2} g_i \zeta_q^i),
\end{displaymath}
which implies that $g_i=g_1\zeta_p^{2\rho}$ for $u^2\times i\equiv 1\modu q$ and so $g_{u_{-2}}=g_1\zeta_p^{2\rho}$.
\item
We continue up to
$\tau^{(q-2)\rho}(g(\mathbf q))=\tau^{q-3}(\zeta_p^\rho g(\mathbf q))=\dots=\zeta_p^{(q-2)\rho} g(\mathbf q)$.
Then
\begin{displaymath}
\sum_{i=1} ^{q-2}g_i \zeta_q^{u^{q-2}i}=\zeta_p^{(q-2)\rho}\times(\sum_{i=1}^{q-2} g_i \zeta_q^i),
\end{displaymath}
which implies that $g_i=g_1\zeta_p^{(q-2)\rho}$ for $u^{q-2}\times i\equiv 1\modu q$
and so $g_{u_{-(q-2)}}=g_1\zeta_p^{(q-2)\rho}$.
\item
Observe that $u$ is a primitive root $\modu q$ and so $u_{-1}$ is a primitive root $\modu q$.
Then it follows that
$g(\mathbf q)
=g_1\times (\zeta_q +\zeta_p^\rho\zeta_q^{u_{-1}}+\zeta_p^{2\rho}\zeta_q^{u_{-2}}+\dots \zeta_p^{(q-2)\rho}\zeta_q^{u_{-(q-2)}})$.
Let $U=\zeta_q +\zeta_p^\rho\zeta_q^{u_{-1}}+\zeta_p^{2\rho}\zeta_q^{u_{-2}}+\dots \zeta_p^{(q-2)\rho}\zeta_q^{u_{-(q-2)}}$.
\item
We prove now that $g_1\in \Z[\zeta_p]^*$. From Stickelberger
relation $g_1^p \times U^p =\mathbf q^{S}$. From
$S=\sum_{i=1}^{p-1}\varpi_t^{-1}\times t$ it follows that
$\varpi_t^{-1}(\mathbf q)^t\ \|\ \mathbf q^{S}$ and so that
$g_1\not\equiv 0\modu \varpi_t^{-1}(\mathbf q)$ because $g_1^p$ is a
$p$-power, which implies that $g_1\in \Z[\zeta_p]^*$. Let us
consider the relation(\ref{e6012211}). Let $x=1\in{\bf F}_q$, then
$Tr_{{\bf F}_q/{\bf F}_q}(x)=1$ and $\chi_\mathbf q^{(p)}(1)=1^{(q-1)/p}\modu \mathbf q$,
so $\chi_\mathbf q^{(p)}(1)=1$ and thus the coefficient
of $\zeta_q$ is $1$ and so $g_1=1$.
\item
From Stickelberger,  $g(\mathbf q)^p \Z[\zeta_p]=\mathbf q^S$,
which achieves the proof.
\end{enumerate}
\end{proof}
\end{lem}
%
\paragraph{Remark:}
$g(\mathbf q)
=\zeta_q +\zeta_p^\rho\zeta_q^{u_{-1}}+\zeta_p^{2\rho}\zeta_q^{u_{-2}}+\dots +\zeta_p^{(q-2)\rho}\zeta_q^{u_{-(q-2)}}$ implies that
$\tau(g(\mathbf q))
=\zeta_q^u +\zeta_p^\rho\zeta_q+\zeta_p^{2\rho}\zeta_q^{u_{-1}}+\dots +\zeta_p^{(q-2)\rho}\zeta_q^{u_{-(q-3)}}$ and that
$\zeta^\rho\times g(\mathbf q)
=\zeta^\rho\zeta_q +\zeta_p^{2\rho}\zeta_q^{u_{-1}}+\zeta_p^{3\rho}\zeta_q^{u_{-2}}+\dots +\zeta_p^{(q-1)\rho}
\zeta_q^{u_{-(q-2)}}$
and we can verify directly that $\tau(g(\mathbf q))=\zeta_p^\rho \times g(\mathbf q)$
for this expression of $g(\mathbf q)$, observing that $q-1\equiv 0\modu p$.
%
When $q\equiv 1\modu p$  Stickelberger's relation is connected with the Kummer's relation on Jacobi resolvents, see for instance
Ribenboim, \cite{rib} (2A) b. p. 118 and (2C) relation (2.6) p. 119.
\begin{lem}\label{l512162}$ $
If $q\equiv 1\modu p$ then
\begin{enumerate}
\item
The value of $g(\mathbf q)$ obtained in relation (\ref{e512151}) is the   Jacobi resolvent:
\be\label{e608051}
g(\mathbf q)=\ <\zeta_p^{-v},\zeta_q>\ =\sum_{i=0}^{q-2}\zeta_p^{-v i}\zeta_q^{u^{-i}}.
\ee
\end{enumerate}
\begin{proof}$ $
\begin{enumerate}
\item
Show first that $g(\mathbf q)=<\zeta_p^{\rho},\zeta_q>$: apply formula of Ribenboim \cite{rib} (2.2) p. 118 with
$p=p,\  q=q,\  \zeta =\zeta_p,\ \rho=\zeta_q, \ n=\rho,\quad u=i,\ m=1$ and $h=u^{-1}$
(where the left members notations $p, q, \zeta,\rho, n, u, m$ and $h$ are the Ribenboim notations).
\item
Show that $\rho=-v$:
we start of
$<\zeta_p^\rho,\zeta_q>=g(\mathbf q)$.
Then $v$ is a primitive root $\modu p$, so there exists  a natural integer $l$ such that
$\rho \equiv v^l\modu p$.
By conjugation $\sigma^{-l}$ we get
$<\zeta_p,\zeta_q>=g(\mathbf q)^{\sigma^{-l}}$.
Raising to $p$-power
$<\zeta_p,\zeta_q>^p=g(\mathbf q)^{p\sigma^{-l}}$.
From lemma \ref{l12161} and Stickelberger relation we get
\be\label{e610011}
<\zeta_p,\zeta_q>^p\Z[\zeta_p]=\mathbf q^{P(\sigma)\sigma^{-l}}.
\ee
From Kummer's relation (2.6) p. 119 in Ribenboim \cite{rib}, we get
$<\zeta_p,\zeta_q>^p\Z[\zeta_p]=\mathbf q^{P_1(\sigma)}$ with $P_1(\sigma)=\sum_{j=0}^{p-2}\sigma^j v_{(p-1)/2-j}$ and thus
\be\label{e610012}
\sum_{i=0}^{p-2} \sigma^{i-l}v_{-i}=\sum_{j=0}^{p-2}\sigma^j v_{(p-1)/2-j}.
\ee
Then, gathering the Kummer and Stickelberger relations (\ref{e610011}) and (\ref{e610012}),  we get
$i-l\equiv j\modu p$,  so $i-j\equiv l\modu p$, in the other hand   $-i\equiv \frac{p-1}{2}-j\modu p$, so $i-j\equiv -\frac{p-1}{2}\modu p$,
so
$l\equiv -\frac{p-1}{2}\modu p$,
so
$l\equiv \frac{p+1}{2}\modu p$,
thus
$\rho\equiv v_{(p+1)/2}\modu p$ and finally $\rho= -v$.
\end{enumerate}
\end{proof}
\end{lem}
%
\begin{lem}\label{l512161}
If $ q\equiv 1\modu p$ then $g(\mathbf q)\equiv -1\modu \pi$.
\begin{proof}
From $g(\mathbf q)=\zeta_q +\zeta_p^{-v}\zeta_q^{u_{-1}}+\zeta_p^{-2v}\zeta_q^{u_{-2}}+\dots +\zeta_p^{-(q-2)v}\zeta_q^{u^{-(q-2)}}$,
we see that
$g(\mathbf q)\equiv \zeta_q +\zeta_q^{u_{-1}}+\zeta_q^{u_{-2}}+\dots +\zeta_q^{u_{-(q-2)}}\modu \pi$.
From $u^{-1}$ primitive root $\modu p$ it follows  that
$1+\zeta_q +\zeta_q^{u_{-1}}+\zeta_q^{u_{-2}}+\dots +\zeta_q^{u_{-(q-2)}}=0$, which leads to the result.
\end{proof}
\end{lem}
%
This result implies that $\pi^p\ |\ g(\mathbf q)^p+1$.
It is possible to improve this    result  by the following  lemma:
\begin{lem}\label{l602061}$ $
If $q\equiv 1\modu p$ and   $p^{(q-1)/p}\not\equiv 1\modu q$ then $\pi^p\ \|\ g(\mathbf q)^p+1$.
\begin{proof}$ $
\begin{enumerate}
\item
We start of $g(\mathbf q)=
\zeta_q +\zeta_p^{-v}\zeta_q^{u_{-1}}+\zeta_p^{-2v}\zeta_q^{u_{-2}}+\dots \zeta_p^{-(q-2)v}\zeta_q^{u_{-(q-2)}}$, so
\begin{displaymath}
g(\mathbf q)=\zeta_q +((\zeta_p^{-v}-1)+1)\zeta_q^{u_{-1}}+((\zeta_p^{-2v}-1)+1)\zeta_q^{u_{-2}},
+\dots ((\zeta_p^{-(q-2)v}-1)+1)\zeta_q^{u_{-(q-2)}}
\end{displaymath}
also
\begin{displaymath}
g(\mathbf q)=-1+(\zeta_p^{-v}-1)\zeta_q^{u_{-1}}+(\zeta_p^{-2v}-1)\zeta_q^{u^{-2}}
+\dots +(\zeta_p^{-(q-2)v}-1)\zeta_q^{u^{-(q-2)}}.
\end{displaymath}
Then $\zeta_p^{-iv}\equiv 1-iv\lambda\modu \pi^2$, so
\begin{displaymath}
g(\mathbf q)\equiv -1+\lambda\times (-v\zeta_q^{u_{-1}}-2v \zeta_q^{u_{-2}}
-\dots -(q-2)v)\zeta_q^{u_{-(q-2)}})\modu\lambda^2.
\end{displaymath}
Then
$g(\mathbf q) =-1+\lambda U+\lambda^2V$
with
\bd
U=-v\zeta_q^{u_{-1}}-2v \zeta_q^{u_{-2}}
+\dots -(q-2)v\zeta_q^{u_{-(q-2)}}
\ed
and $U,V\in\Z[\zeta_{pq}]$.
\item
Suppose that $\pi^{p+1} \ |\ g(\mathbf q)^p+1$ and search for a contradiction:
then, from $g(\mathbf q)^p =(-1+\lambda U+\lambda^2V)^p$,
it follows that
 $p\lambda U+\lambda^pU^p\equiv 0\modu\pi^{p+1}$
and so $U^p-U\equiv 0\modu \pi$ because $p\lambda+\lambda^p\equiv 0\modu\pi^{p+1}$.
Therefore
\begin{displaymath}
\begin{split}
& (-v\zeta_q^{u_{-1}}-2v \zeta_q^{u_{-2}}
-\dots -(q-2)v\zeta_q^{u_{-(q-2)}})^p-\\
&(-v\zeta_q^{u_{-1}}-2v \zeta_q^{u_{-2}}
+\dots -(q-2)v\zeta_q^{u_{-(q-2)}})\equiv 0\modu \lambda,\\
\end{split}
\end{displaymath}
and so
\begin{displaymath}
\begin{split}
& (\zeta_q^{pu_{-1}}+2 \zeta_q^{pu_{-2}}
+\dots +(q-2)\zeta_q^{pu_{-(q-2)}})\\
& -(\zeta_q^{u_{-1}}+2 \zeta_q^{u_{-2}}
+\dots +(q-2)\zeta_q^{u_{-(q-2)}})\equiv 0\modu \lambda.\\
\end{split}
\end{displaymath}
\item
For any  natural  $j$ with $1\leq j\leq q-2$,  there must exist a natural $j^\prime$ with $1\leq j^\prime\leq q-2$ such that simultaneously:
\begin{displaymath}
\begin{split}
& pu^{-j^\prime}\equiv u^{-j}\modu q\Rightarrow p\equiv u^{j^\prime-j}\modu q,\\
& j^\prime \equiv   j\modu p.\\
\end{split}
\end{displaymath}
Therefore $p\equiv u^{p\times \{(j^\prime-j)/p\}}\modu q$
and so $p^{(q-1)/p}\equiv u^{p\times (q-1)/p)\times \{(j^\prime-j)/p\}}\modu q$
thus $p^{(q-1)/p}\equiv 1\modu q$, contradiction.
\end{enumerate}
\end{proof}
\end{lem}
%
Let $M$ be the subfield of $K_{pq}$ with degree $[M:K]=p$. We know that   $g(\mathbf q)^p\in K$ and $g(\mathbf q)\not\in K$, thus
$M= K(g(\mathbf q))$ is a  cyclic extension of degree $p$.
In the following lemma we connect the decomposition of the prime ideal $\pi$ of $K$  in $M$ with the value of $p^{(q-1)/p}\modu q$
used in  lemma \ref{l602061} p. \pageref{l602061}.
\begin{lem}\label{l612011}$ $
\bn
\item
$\pi$ does not ramify in $M/K$.
\item
If $\pi$ is inert  in $M/K$ then $\pi^p\ \|\ g(\mathbf q)^p+1$ and $p^{(q-1)/p}\not\equiv 1\modu p$.
\item
If $\pi$  splits  in $M/K$ then $\pi^{p+1}\ |\ g(\mathbf q)^p+1$ and  $p^{(q-1)/p}\equiv 1\modu p$.
\en
\begin{proof}$ $
\bn
\item
$g(\mathbf q)^p+1\equiv 0\modu\pi^p$, therefore  $g(\mathbf q)$ is primary and $\pi$ does not ramify in $M/K$ (see for instance
Washington \cite{was} exercise 9.3 (b) p. 183).
\item
We start of
$g(\mathbf q)=\sum_{i=0}^{q-2}\zeta_p^{-vi}\zeta_q^{u_{-i}}\in M$.
Let $\tau$ be the $K$-isomorphism of $K_{pq}$ given by $\tau :\zeta_q\rightarrow\zeta_q^{u_{-1}}$.
Let us note  $\sigma$ the extension of the $\Q$-isomorphism $\sigma$ of $K$ to $M$ given by $\sigma :\zeta_q\rightarrow \zeta_q$.
\item
Let $G=\sigma^{-1+(p-1)/2}\circ\tau^{-1}(g(\mathbf q))$.
Then $G=\sum_{i=0}^{q-2}\zeta_p^{i}\zeta_q^{u_i}\in M$.
Without loss of generality, we work now with $G$ in place of $g(\mathbf q)$.
Thus
$G=\sum_{i=0}^{q-2}(\lambda+1)^i\zeta_q^{u_i}$,
hence
$G\equiv \sum_{i=0}^{q-2}\zeta_q^{u_i}+\lambda\sum_{i=0}^{q-2} i\zeta_q^{u_i}\modu\pi^2$,
hence
$G\equiv -1+\lambda\sum_{i=0}^{q-2} i\zeta_q^{u_i}\modu\pi^2$.
\item
Let $\Pi$ be a prime of $M$ lying over $\pi$.
$\frac{G+1}{\lambda}\equiv \sum_{i=0}^{q-2}i\zeta_q^{u_i}\modu\Pi$.
\bn
\item
\underline{If $\pi$ does not split in $M/K$} then  $\tau(\Pi)=\Pi$.
If
$\frac{G+1}{\lambda}\equiv 0\modu\Pi$ then  $\frac{G\zeta+1}{\lambda}\equiv 0\modu\tau(\Pi)$,
thus   $\frac{G\zeta+1}{\lambda}\equiv 0\modu\Pi$
hence $G\equiv 0\modu\Pi$,  contradiction, hence $\frac{G+1}{\lambda}$ is coprime with $p$ and
$\pi^p\ \|\ G^p+1$.
\item
\underline{If $\pi$ splits in $M/K$}  then there exists a natural $w$ such that $g(\mathbf q)+\zeta^w\equiv 0\modu \Pi^2$
where $\Pi$ is a prime of $M$ lying over $\pi$ because $\Pi$ is of inertial degree $1$ and $\pi^{p+1}\ |\ G^p+1$.
\en
\item
Then apply lemma \ref{l602061} p.\pageref{l602061}.
\en
\end{proof}
\end{lem}
%
%
\subsection{The polynomial   $P(\sigma)=\sum_{i=0}^{p-2}\sigma^i v_{-i}$.}
This subsection contains some elementary algebraic computations on $P(\sigma)\in\Z[G_p]$.
%
\begin{lem}\label{l512165}
\begin{equation}\label{e512172}
P(\sigma)\times (\sigma-v)= p\times Q(\sigma),
\end{equation}
where $Q(\sigma)=\sum_{i=1}^{p-2}\delta_i\times \sigma^i\in\Z[G_p]$ is given by
\begin{equation}
\begin{split}
& \delta_{p-2}= \frac{v_{-(p-3)}-v_{-(p-2)}v}{p},\\
& \delta_{p-3}= \frac{v_{-(p-4)}-v_{-(p-3)}v}{p},\\
& \vdots\\
& \delta_i=\frac{v_{-(i-1)}-v_{-i} v}{p},\\
&\vdots\\
& \delta_1 = \frac{1-v_{-1}v}{p},\\
\end{split}
\end{equation}
with $-p<\delta_i\leq 0$.
\begin{proof}
We start of  the relation
$P(\sigma)= \prod_{i=0,\ i\not=1}^{p-2} (\sigma-v_i)$ in ${\bf F}_p[G_p]$, hence
$P(\sigma)\times(\sigma-v)=p\times Q(\sigma)$,
with $Q(\sigma)\in\Z[G_p]$.
Then  we identify in $\Z[G_p]$ the  coefficients in the relation
\begin{displaymath}
\begin{split}
&(v^{-(p-2)}\sigma^{p-2}+v^{-(p-3)}\sigma^{p-3}+\dots+v^{-1}\sigma+1)\times(\sigma-v)=\\
&p\times (\delta_{p-2}\sigma^{p-2}+\delta_{p-3}\sigma^{p-3}+\dots+\delta_1\sigma+\delta_0),\\
\end{split}
\end{displaymath}
where $\sigma^{p-1}=1$.
\end{proof}
\end{lem}
\paragraph{Remark:}
\begin{enumerate}
\item
Observe that, with our notations,  $\delta_i\in \Z,\quad i=1,\dots,p-2$, but generally $\delta_i\not\equiv 0\modu p$.
\item
We see also that $-p< \delta_i\leq 0$.
Observe also that $\delta_0=\frac{v^{-(p-2)}-v}{p}=0$.
\end{enumerate}
%
\begin{lem}\label{l601211}
The polynomial $Q(\sigma)$ verifies
\begin{equation}\label{e6012110}
Q(\sigma)=\{(1-\sigma)(\sum_{i=0}^{(p-3)/2}\delta_i\times \sigma^i)+(1-v)\sigma^{(p-1)/2}\}\times(\sum_{i=0}^{(p-3)/2}\sigma^i).
\end{equation}
\begin{proof}
We start of $\delta_i=\frac{v^{-(i-1)}-v^{-i}v}{p}$.
Then
\begin{displaymath}
\delta_{i+(p-1)/2}=\frac{v^{-(i+(p-1)/2-1)}-v^{-(i+(p-1)/2)}v}{p}=
\frac{p-v^{-(i-1)}-(p-v^{-i})v}{p}=1-v-\delta_i.
\end{displaymath}
Then
\begin{displaymath}
\begin{split}
& Q(\sigma)=\sum_{i=0}^{(p-3)/2} (\delta_i\times (\sigma^i-\sigma^{i+(p-1)/2}+(1-v)\sigma^{i+(p-1)/2})\\
& =(\sum_{i=0}^{(p-3)/2}\delta_i\times\sigma^i)\times (1-\sigma^{(p-1)/2)})
+(1-v)\times \sigma^{(p-1)/2}\times (\sum_{i=0}^{(p-3)/2}\sigma^i),\\
\end{split}
\end{displaymath}
which leads to the result.
\end{proof}
\end{lem}
Let us note in the sequel
\begin{equation}\label{e608191}
Q_1(\sigma)=(1-\sigma)(\sum_{i=0}^{(p-3)/2}\delta_i\times \sigma^i)+(1-v)\sigma^{(p-1)/2}.
\end{equation}
%
\subsection{$\pi$-adic congruences on the singular integers $A$}
The results of this section are algebraic preliminaries to the two following sections \ref{s610011} p. \pageref{s610011}
and \ref{s610012} p. \pageref{s610012}.

We remind that the group of ideal classes of the cyclotomic field
is generated by the ideal classes of prime ideals of degree $1$, see for instance Ribenboim, \cite{rib} (3A) p. 119.
The prime ideal $\mathbf q$ of $\Z[\zeta_p]$ has a non-trivial class $Cl(\mathbf q)\in \Gamma$ where $\Gamma$ is
a subgroup of order $p$ of $C_p$  defined in section \ref{s601191} p. \pageref{s601191},
with a singular integer $A$  given  by $A \Z[\zeta_p] = \mathbf q^p$.
%
\begin{lem}\label{l512163}$ $
One has $g(\mathbf q)^{2p^2}\eta=A^{P(\sigma)}$ and
$(\frac{g(\mathbf q}{\overline{g(\mathbf q})})^{2p^2}=(\frac{A}{\overline{A}})^{P(\sigma)}$
with some unit $\eta\in\Z[\zeta_p+\zeta_p^{-1}]^*$.
\begin{proof}
We start of
$g(\mathbf q)^p \Z[\zeta_p]= \mathbf q^{P(\sigma)}$.
Raising to $p$-power we get
$g(\mathbf q)^{2p^2} \Z[\zeta_p]= \mathbf q^{2pP(\sigma)}$.
But $A \Z[\zeta_p] = \mathbf q^{2p}$, so
\begin{equation}\label{e601071}
g(\mathbf q)^{2p^2} \Z[\zeta_p]= A^{P(\sigma)}\Z[\zeta_p],
\end{equation}
so
\begin{equation}\label{e601063}
 g(\mathbf q)^{2p^2}\times  \zeta_p^w\times \eta= A^{P(\sigma)},\quad \eta\in \Z[\zeta_p+\zeta_p^{-1}]^*,
\end{equation} where $w$ is a natural number.
Therefore, by complex conjugation, we get
$ \overline{g(\mathbf q)}^{2p^2}\times\zeta_p^{-w}\times  \eta= \overline{A}^{P(\sigma)}.$
Then
$ (\frac{g(\mathbf q)}{\overline{g(\mathbf q)}})^{p^2}\times\zeta_p^{2w}=(\frac{A}{\overline{A}})^{P(\sigma)}$.
$A$ is semi-primary, thus  $w=0$.
Then $ (\frac{g(\mathbf q)}{\overline{g(\mathbf q)}})^{2p^2}=(\frac{A}{\overline{A}})^{P(\sigma)}$.
\end{proof}
\end{lem}
\paragraph{Remark:} Observe that this lemma is true if either $q\equiv 1\modu p$ or $q\not\equiv 1\modu p$.
%
\begin{thm}\label{t601311}$ $
\begin{enumerate}
\item
$g(\mathbf q)^{2p^2}=\pm A^{P(\sigma)}$.
\item
$g(\mathbf q)^{2p(\sigma-1)(\sigma-v)}=\pm
(\frac{\overline{A}}{A})^{Q_1(\sigma)}$
where
\begin{displaymath}
Q_1(\sigma)= (1-\sigma)\times (\sum_{i=0}^{(p-3)/2}\delta_i\times \sigma^i)+(1-v)\times \sigma^{(p-1)/2}.
\end{displaymath}
\end{enumerate}
\begin{proof}$ $
\begin{enumerate}
\item
We start of  $g(\mathbf q)^{2p^2}\times\eta = A^{P(\sigma)}$ proved.
Then $g(\mathbf q)^{2p^2(\sigma-1)(\sigma-v)}\times\eta^{(\sigma-1)(\sigma-v)}=
A^{P(\sigma)(\sigma-1)(\sigma-v)}$.
From lemma \ref{l601211}, we get
\begin{displaymath}
P(\sigma)\times (\sigma-v)\times (\sigma-1)=p \times Q_1(\sigma)\times (\sigma^{(p-1)/2}-1).
\end{displaymath}
Therefore
\begin{equation}\label{e602013}
g(\mathbf q)^{2p^2(\sigma-1)(\sigma-v)}\times\eta^{(\sigma-1)(\sigma-v)}=
(\frac{\overline{A}}{A})^{p Q_1(\sigma)},
\end{equation}
and by conjugation
\begin{displaymath}
\overline{g(\mathbf q)}^{2p^2(\sigma-1)(\sigma-v)}\times\eta^{(\sigma-1)(\sigma-v)}=
(\frac{A}{\overline{A}})^{p Q_1(\sigma)}.
\end{displaymath}
Multiplying these two relations we get, observing that $g(\mathbf q)\times \overline{g(\mathbf q)}=q^f$,
\begin{displaymath}
q^{2f p^2(\sigma-1)(\sigma-v)}\times \eta^{2(\sigma-1)(\sigma-v)}=1,
\end{displaymath}
also
\begin{equation}\label{e611041}
\eta^{2(\sigma-1)(\sigma-v)}=1.
\end{equation}
Show that it implies $\eta=\pm 1$:
If $\eta^{\sigma-1}\not \in\Q$ then $v_\pi(\eta^{\sigma-1}-1)=2n$ for a  natural $n\geq 1$.
This implies that $v_\pi(\eta^{(\sigma-1)(\sigma-v)}-1)=v_\pi(\eta^{\sigma-1}-1)=2n$ because the action of  $\sigma-v$ cannot change
the parity of a number of even $\pi$-adic parity, hence that $\eta^{(\sigma-1)(\sigma-v)}\not=1$, contradiction.
Therefore $\eta\in\Q$ and relation (\ref{e611041}) follows.
This leads, with relation (\ref{e602013}),  to $g(\mathbf q)^{2p^2}=\pm A^{P(\sigma)}$ and achieves the proof of the first part.
\item
From relation (\ref{e602013}) we get
\begin{equation}\label{e602014}
g(\mathbf q)^{2p^2(\sigma-1)(\sigma-v)}=\pm
(\frac{\overline{A}}{A})^{p Q_1(\sigma)},
\end{equation}
so
\begin{equation}\label{e602031}
g(\mathbf q)^{2p(\sigma-1)(\sigma-v)}=\pm \zeta_p^w\times
(\frac{\overline{A}}{A})^{ Q_1(\sigma)},
\end{equation}
where $w$ is a natural number.
But $g(\mathbf q)^{\sigma-v}\in K_p$ and so $g(\mathbf q)^{p(\sigma-v)(\sigma-1)}\in (K_p)^p$,
see for instance Ribenboim \cite{rib} (2A) b. p. 118.
and $(\frac{\overline{A}}{A})^{Q_1(\sigma)}\in (K_p)^p$ because
$\sigma-\mu \ |\ Q_1(\sigma)$ in ${\bf F}_p[G_p]$ when $A$ is singular negative
and $\frac{\overline{A}}{A}=D^p$ when $A$ is singular positive imply that $w=0$,
which achieves the proof of the second part.
\end{enumerate}
\end{proof}
\end{thm}
%
\paragraph{Remarks}
\begin{enumerate}
\item
Observe that this theorem is true either $q\equiv 1\modu p$ or $q\not\equiv 1\modu p$.
\item
$g(\mathbf q)\equiv -1\modu\pi$ implies that
$g(\mathbf q)^{p^2}\equiv -1\modu \pi$.
Observe that if $A\equiv a\modu\pi$ with $a$ natural number then
$A^{P(\sigma)}\equiv a^{1+v^{-1}+\dots+v^{-(p-2)}}=a^{p(p-1)/2}\modu \pi\equiv \pm 1\modu \pi$ consistent with previous result.
\end{enumerate}
%
\begin{thm}\label{l512164}$ $

\begin{enumerate}
\item
If $q\equiv 1 \modu p$ then $ A^{P(\sigma)}\equiv \delta \modu \pi^{2p-1}$ with $\delta\in\{-1,1\}$.
\item
If $q\not\equiv 1 \modu p$ then $ A^{P(\sigma)}\equiv \delta\modu \pi^{2p}$ with $\delta\in\{-1,1\}$.
\item
If and only if $q\equiv 1 \modu p$ and $p^{(q-1)/p}\not\equiv 1\modu q$ then $\pi^{2p-1}\ \|\   A^{P(\sigma)} -\delta $ with $\delta\in\{-1,1\}$.
\end{enumerate}
\begin{proof}$ $
\begin{enumerate}
\item
From lemma \ref{l512161}, we get $\pi^p\ |\ g(\mathbf q)^p+1$ and so $\pi^{2p-1}\ |\ g(\mathbf q)^{p^2}+1$.
Then apply theorem \ref{t601311}.
\item
From lemma \ref{l512151}, then $g(\mathbf q)\in\Z[\zeta_p]$ and so $\pi^{p+1}\ |\ g(\mathbf q)^p+1$ and also
$\pi^{2p}\ |\ g(\mathbf q)^{p^2}+1$.
\item
Applying  lemma \ref{l602061} we get $\pi^p\ \|\ g(\mathbf q)^p+1$ and so $\pi^{2p-1}\ \|\ g(\mathbf q)^{p^2}+1$. Then
apply theorem \ref{t601311}.
\end{enumerate}
\end{proof}
\end{thm}
%
\paragraph{Remark:} If $C\in\Z[\zeta_p]$ is any semi-primary number with $C\equiv c\modu \pi^2$ we can only assert in general that
$v_\pi(C^{P(\sigma)}-c_{p-1})\geq p-1$ for $c_{p-1}\in\Z$: we start from $P(\sigma)=\prod_{i=0,\ i\not=1}^{p-2} (\sigma-v_i)+p R(\sigma)$ where
$R(\sigma)\in\Z[\zeta_p]$. Then $v_\pi(C-c)=m$ with  $1<m\leq p-2$, hence $v_\pi(C^{\sigma-v_m}-c_m)=m_1>m$. Pursuing to climb up,
$v_\pi(C^{(\sigma-v_m)(\sigma-v_{m_1})}-c_{m_1})=m_2>m_1,\dots$,
 we obtain the result.
For the singular numbers $A$ considered, we assert more $v_{\pi}(A^{P(\sigma)}\pm 1)\geq 2p-1$. We obtain a $\pi$-adic depth
$\geq 2p-1$ for singular numbers $A$ where the $\pi$-adic depth is only $\geq p-1$ in the general case of semi-primary numbers $C$.
%
\section{Jacobi resolvents and Jacobi cyclotomic function}\label{s610011}
The $\pi$-adic development of singular numbers shall be described in two steps: application of Jacobi resolvents and
then of Jacobi cyclotomic functions.
%
\subsection{Singular  integers  and Jacobi resolvents  }\label{s608221}
\bn
\item
Let $\mathbf q$ be a prime ideal of $\Z[\zeta_p]$ of inertial degree $1$. Let $q$ be the prime number lying above $\mathbf q$ so with
$q\equiv 1\modu p$.
We can suppose, from Kummer,  that $Cl(\mathbf q)\in\Gamma\in C_p$ where the group $\Gamma$ of order $p$ is annihilated
by $\sigma-\mu$ with $\mu=v_n\in {\bf F}_p^*$ where $1<n\leq p-2$ because $\mu\not=v$ from Stickelberger theorem.
Let $A$ be a singular  number with $A\Z[\zeta_p]=\mathbf q^{2p}$.
We recall that $A^{\sigma-\mu}=\alpha^p$ for $\alpha\in K_p$ and $v_\pi(\frac{A}{\overline{A}}-1)\geq 3$.

\item
$v_{\pi}((\frac{\alpha}{\overline{\alpha}})^p-1)\geq p-1+3$.
Let us define the  natural integer $\nu\geq 2$ by
\be\label{e609203}
\begin{split}
& \nu\equiv  v_{\pi}((\frac{\alpha}{\overline{\alpha}})^p-1)\modu (p-1), \ 3\leq \nu\leq p-2.\\
\end{split}
\ee
\item
$g(\mathbf q)\not\in K_p$, but $g(\mathbf q)^{\sigma-v}\in K_p$, (see for instance Ribenboim \cite{ri2}, F. p. 440.
Then let us define the natural integer  $\rho$ by
\be\label{e608291}
\begin{split}
& \rho =v_\pi((\frac{g(\mathbf q)}{\overline{g(\mathbf q)}})^{p(\sigma-v)}-1)-(p-1).\\
\end{split}
\ee
The definition relation (\ref{e608291})  of $\rho$
implies that $\rho$ is odd and that $\rho\geq 3$.
\en
%
\clearpage
The following theorem connects the  $\pi$-adic expansion of singular integers $A$ with Jacobi resolvents $g(\mathbf q)$:
\begin{thm}\label{t608201}$ $
\bn
\item
If $n\not=\rho$ then $\nu\leq\rho$.
\item
If $A$ is singular primary then
$v_\pi(\frac{A}{\overline{A}}-1)\leq p-1+\rho$.
\en
\begin{proof}$ $
\bn
\item
From lemma  \ref{l512163} p. \pageref{l512163}
we obtain $(\frac{g(\mathbf q)}{\overline{g(\mathbf q)}})^{2p^2}=(\frac{A}{\overline{A}})^{P(\sigma)}$, so
$(\frac{g(\mathbf q)}{\overline{g(\mathbf q)}})^{2p^2(\sigma-v)}=(\frac{A}{\overline{A}})^{P(\sigma)(\sigma-v)}$.
From lemma \ref{l512165} p. \pageref{l512165} we obtain
$(\frac{g(\mathbf q)}{\overline{g(\mathbf q)}})^{2p^2(\sigma-v)}=(\frac{A}{\overline{A}})^{p Q(\sigma)}$,
so
\bd
(\frac{g(\mathbf q)}{\overline{g(\mathbf q)}})^{2p(\sigma-v)}=(\frac{A}{\overline{A}})^{ Q(\sigma)},
\ed
also
\bd
(\frac{g(\mathbf q)}{\overline{g(\mathbf q)}})^{2p(\sigma-v)(1-\sigma)}=(\frac{A}{\overline{A}})^{ Q(\sigma)(1-\sigma)}.
\ed
From lemma \ref{l601211} p. \pageref{l601211} and relation (\ref{e608191}) p. \pageref{e608191}  we get
\bd
(\frac{g(\mathbf q)}{\overline{g(\mathbf q)}})^{2p(\sigma-v)(1-\sigma)}=(\frac{A}{\overline{A}})^{ Q_1(\sigma)(1-\sigma^{(p-1)/2})},
\ed
so
\be\label{e608215}
(\frac{g(\mathbf q)}{\overline{g(\mathbf q)}})^{2p(\sigma-v)(1-\sigma)}=(\frac{A}{\overline{A}})^{2 Q_1(\sigma))},
\ee
and
\be\label{e608292}
(\frac{g(\mathbf q)}{\overline{g(\mathbf q)}})^{2p(\sigma-v)(1-\sigma)(\sigma-\mu)}=(\frac{\alpha}{\overline{\alpha}})^{2 p Q_1(\sigma))}.
\ee
\item
From lemma \ref{l512161} p. \pageref{l512161} we know that $\frac{g(\mathbf q)}{\overline{g(\mathbf q)}}\equiv 1\modu\pi$.
From Ribenboim \cite{ri2} F. p. 440 $(\frac{g(\mathbf q)}{\overline{g(\mathbf q)}})^{\sigma-v}\in K_p$.
Therefore
\bd
(\frac{g(\mathbf q)}{\overline{g(\mathbf q)}})^{p(\sigma-v)}\equiv 1\modu\pi^{p-1+2}.
\ed
From the definition of $\rho$ and the hypothesis $n\not=\rho$ it follows that
\be\label{e608231}
v_\pi( (\frac{g(\mathbf q)}{\overline{g(\mathbf q)}})^{p(\sigma-v)(1-\sigma)(\sigma-\mu)}-1)= p-1+\rho,
\ee
because $v_{\pi}((\frac{g(\mathbf q)}{\overline{g(\mathbf q)}})^{p(\sigma-v)}=p-1+\rho$ with  $\rho$ odd and $\rho\not= n$ implies
$\mu\not= v_\rho$;
indeed $(\frac{g(\mathbf q)}{\overline{g(\mathbf q)}})^{p(\sigma-v)}$ cannot climb up by the action of $(\sigma-1)(\sigma-\mu)$.
Thus $v_\pi((\frac{\alpha}{\overline{\alpha}})^p-1)\leq p-1+\rho$ follows from relation (\ref{e608292}).
\item
$v_\pi(\frac{A}{\overline{A}}-1)\leq p-1+\rho$ is an immediate consequence of relation (\ref{e608215}).
\en
\end{proof}
\end{thm}
%
%
\paragraph{Remark:} Observe that $g(\mathbf q)$ is explicitly computable by $g(\mathbf q)=\sum_{i=0}^{p-2}\zeta_p^{-vi}\zeta_q^{u^{-i}}$ and that
the number $\rho$ defined in relation (\ref{e608291}) p. \pageref{e608291} can therefore be explicitly be algebraically computed by
\bd
\rho=v_\pi((\frac{\sum_{i=0}^{p-2}\zeta_p^{-vi}\zeta_q^{u^{-i}}}{\sum_{i=0}^{p-2}\zeta_p^{vi}\zeta_q^{-u^{-i}}})^{\sigma-v}-1).
\ed
Numerical MAPLE computations for a large number of pairs of odd prime $(p,q), \ q\equiv 1\modu p$ had shown that $\rho=3$ for almost all pairs
and that $\rho$ seems always small before $p$ when $p$ is large.
%

%
%
\subsection{Jacobi cyclotomic function and singular numbers }\label{s609201}
The  computation method applied is  derived  from Jacobi cyclotomic function (see for instance  Ribenboim \cite{rib} section 2. p. 117).
Recall that $\nu$ and $\rho$ are respectively defined in relations (\ref{e609203}) and (\ref{e608291}).
We shall use the inequality $\nu\leq\rho$ found above
%

\paragraph{Definition:}
let $i$ be an integer with $1\leq i\leq q-2$. There exists one and only one integer $s,\ 1\leq s\leq q-2$ such that
$i\equiv u^s\modu q$. The number s is called the {\it index} of $i$ {\it relative} to $u$ and denoted $s=ind_u(i)$.
Let $a,b$ be natural  numbers  such that $ab(a+b)\not\equiv 0\modu p$.
The Jacobi resolvents verify the relation:
\begin{equation}\label{e605031}
\frac{<\zeta_p^a,\zeta_q><\zeta_p^b,\zeta_q>}{<\zeta_p^{a+b},\zeta_q>}
=
\sum_{i=1}^{q-2}\zeta_p^{a\times  ind_u(i)-(a+b)\times ind_u(i+1)},
\end{equation}
(See for this result for instance Ribenboim \cite{ri2} proposition (I) p. 442).
The interest of this formula for $\pi$-adic structure of singular primary numbers is that the right member belongs to $\Z[\zeta_p]$
though  the Jacobi resolvents $<\zeta_p^a,\zeta_q>,\ <\zeta_p^b,\zeta_q>$ and $<\zeta_p^{a+b},\zeta_q>$ are only in $\Z[\zeta_p,\zeta_q]$:
the left member in $\Z[\zeta_p]$ is more convenient for studying $\pi$-adic congruences.
%
\begin{thm}\label{t605031}
Let $a,b$ be two natural numbers such that $ab\times (a+b)\not\equiv 0\modu p$.
Let $\mathbf q$ be a prime ideal of $\Z[\zeta_p]$ of inertial degree $1$
with $Cl(\mathbf q)\in C_p$ annihilated by $\sigma-\mu,\ \mu\in{\bf F}_p^*$.
Let $q$ be the prime number above $\mathbf q$.
Let $A$ be a corresponding singular number with $A\Z[\zeta_p]=\mathbf q^{2p}$.
If $\rho\leq p$ then
\begin{enumerate}
\item
The Jacobi cyclotomic function
\begin{equation}\label{e605032}
\psi_{a,b}(\zeta_p)=\sum_{i=1}^{q-2}\zeta_p^{a\times ind_u(i)-(a+b)\times ind_u(i+1)}\equiv -1\modu\pi^\nu.
\end{equation}
\item
For $k=2,\dots,\nu-1$
\begin{equation}\label{e605033}
\sum_{i=1}^{q-2} \{a\times ind_u(i)-(a+b)\times ind_u(i+1)\}^k\equiv 0\modu p.
\end{equation}
\end{enumerate}
\begin{proof}$ $
\begin{enumerate}
\item
We start of $g(\mathbf q)=\sum_{i=1}^{q-2}\zeta_p^{-iv}\zeta_q^{u^{-i}}$. There exists  a natural number  $\alpha$ such that
$-v^{\alpha+1}\equiv a\modu p$.
Then $<\zeta_p^a,\zeta_q>= <\zeta_p^{-v v^\alpha},\zeta_q>=\sigma^\alpha(g(\mathbf q)).$
From definition relation (\ref{e608291} p. \pageref{e608291}, $g(\mathbf q)^{\sigma^\alpha-v^\alpha}\equiv \pm 1\modu \pi^\rho$.
Similarly, for $-v^{\beta+1}\equiv b\modu p$,  $<\zeta_p^b,\zeta_q>= \sigma^\beta(g(\mathbf q))$ and
$g(\mathbf q)^{\sigma^\beta-v^\beta}\equiv \pm 1\modu \pi^\rho$.
\item
$<\zeta_p^{a+b},\zeta_q>=<\zeta_p^{-vv^\gamma},\zeta_q>$ with
$-v^{\gamma+1}\equiv a+b\modu p$.
Then
$<\zeta_p^{a+b},\zeta_q>=\sigma^\gamma(g(\mathbf q))$
with
$g(\mathbf q)^{\sigma^\gamma-v^\gamma}\equiv \pm 1\modu \pi^\rho$.
\item
\begin{displaymath}
\frac{<\zeta_p^a,\zeta_q><\zeta_p^b,\zeta_q>}{<\zeta_p^{a+b},\zeta_q>}
=\frac{\sigma^\alpha(g(\mathbf q))\sigma^\beta(q(\mathbf q))}{\sigma^\gamma(g(\mathbf q))}
\equiv \pm g(\mathbf q)^{v^\alpha+v^\beta-v^\gamma}\modu\pi^\rho,
\end{displaymath}
also
\begin{displaymath}
\frac{<\zeta_p^a,\zeta_q><\zeta_p^b,\zeta_q>}{<\zeta_p^{a+b},\zeta_q>}
\equiv \pm g(\mathbf q)^{v^{-1}\times (v^{\alpha+1}+v^{\beta+1}-v^{\gamma+1})\modu p}\modu\pi^\rho.
\end{displaymath}
But
$v^{\alpha+1}+v^{\beta+1}-v^{\gamma+1}\equiv -a-b+a+b\equiv 0\modu p$.
From $g(\mathbf q)\equiv -1\modu p$ it follows that $g(\mathbf q)^p\equiv -1\modu\pi^p$.
Then if $\rho\leq p$ we get
\begin{displaymath}
\frac{<\zeta_p^a,\zeta_q><\zeta_p^b,\zeta_q>}{<\zeta_p^{a+b},\zeta_q>}
\equiv \pm 1\modu\pi^\rho,
\end{displaymath}
and from Jacobi cyclotomic function relation (\ref{e605031})
\begin{displaymath}
\sum_{i=1}^{q-2}\zeta_p^{a\times ind_u(i)-(a+b)\times ind_u(i+1)}\equiv \pm 1\modu\pi^\rho.
\end{displaymath}
But we see directly , from $\zeta_p\equiv 1\modu p$ and $q-2\equiv -1\modu p$, that
\begin{displaymath}
\sum_{i=1}^{q-2}\zeta_p^{a\times ind_u(i)-(a+b)\times ind_u(i+1)}\equiv -1\modu\pi,
\end{displaymath}
and finally that
\begin{displaymath}
\sum_{i=1}^{q-2}\zeta_p^{a\times ind_u(i)-(a+b)\times ind_u(i+1)}\equiv -1\modu\pi^\rho,
\end{displaymath} which proves relation (\ref{e605032}) because  $\nu\leq \rho$ from  theorem
\ref{t608201}.
\item
The congruences (\ref{e605033}) are an immediate consequence,  using logarithmic derivatives.
\end{enumerate}
\end{proof}
\end{thm}
%
This result takes a very simple form when $a=1$ and $b=-2$.
\begin{cor}\label{t605061}
With $a=1,\ b=-2$ in theorem \ref{t605031} p. \pageref{t605031}, if $\rho\leq p$ then
\begin{equation}\label{e605061}
\psi_{1,-2}(\zeta_p)=\sum_{i=1}^{q-2}\zeta_p^{ind_u(i(i+1))}\equiv -1\modu\pi^\nu.
\end{equation}
For $k=2,\dots,\nu-1$
\begin{equation}\label{e609211}
\sum_{i=1}^{q-2} ind_u(i(i+1))^k\equiv 0\modu p.
\end{equation}
\begin{proof}$ $
Take $a=1$ and $b=-2$ in congruence  (\ref{e605032}) to get
\begin{equation}\label{e605062}
\sum_{i=1}^{q-2}\zeta_p^{ind_u(i)+ind_u(i+1)}\equiv -1\modu\pi^\nu.
\end{equation}
But  $ind_u(i)+ind_u(i+1)\equiv ind_u(i(i+1))\modu (q-1)$ and so from $q\equiv 1\modu p$
we get $ind_u(i)+ind_u(i+1)\equiv ind_u(i(i+1))\modu p$ which achieves the proof.
\end{proof}
\end{cor}
\paragraph{Remarks:}
\bn
\item
These results are meaningful for  the case of  singular  primary  and the case of singular  non-primary integers $A$.
\item
These results are meaningful for all prime ideals $\mathbf q$ of inertial degree $1$ whose class is in $\Gamma$.
Observe there are infinitely many such ideals $\mathbf q$.
\item
We have computed $\psi_{1,-2}(\zeta_p)$ for a large number of pairs $(p,q)$ with small $q,\ p\not= 3,\ q\not\equiv 1\modu 3$
and we have found that for almost all  these pairs
$\pi^3\ \|\  \psi_{1,-2}(\zeta_p)$. We conjecture that $\rho$ is in most cases  {\it small} before $p$. Thus we conjecture that
$\nu$ is in most cases {\it small} before $p$.
\en
%
\paragraph{A question:} From Denes it is known that there exists a fundamental system of units $\varepsilon_i\in O_{K^+}^*$
with $\varepsilon_i\equiv  a_i+b_i\lambda^{2i+c_i(p-1)}\modu \pi^{2i+1+c_i(p-1)}$ for $i=1,\dots,\frac{p-3}{2}$
where $a_i,b_i,c_i$ are natural integers and $c_i>0$ when
$\varepsilon_i$ is a primary unit (see for instance Ribenboim \cite{rib} (8D) p. 192).
Is such a result can be extended to singular primary numbers?
In that case it should mean that $\rho$ defined by relation (\ref{e608291})  p. \pageref{e608291} could take
a value not always small in regards of $p$.
By opposite,  our MAPLE computations of $\rho$ seems to show that $\rho=3$ for almost prime numbers $q\equiv 1\modu p$
and that $\rho$ seems to be always small in regards to $p$.
If the Denes result could be extended to all singular numbers $A$ non-units it could explain the scarcity (or impossibility) to
find components of the positive $p$-class group $C_p^+$, (a question  connected to Vandiver's conjecture).

%
\section{Principal prime ideals of $K_p$ and Kummer-Stickelberger relation}\label{s610012}

Observe that the Stickelberger relation and its consequences on prime ideals $\mathbf q$ of $\Z[\zeta_p]$  is meaningful even if
 $\mathbf q$ is a principal ideal.
 As an example, we give the  theorem:
%
\begin{thm}\label{t602061}
Let $q_1\in\Z[\zeta_p]$ be an  integer verifying  $q_1\equiv a\modu \pi^{p+1}$ where $a\in\Z,\quad  a\not\equiv 0\modu p$.
If   $q=N_{K/\Q}(q_1)$ is a prime number then  $p^{(q-1)/p}\equiv 1\modu q$.
\begin{proof}
$a^{p-1}\equiv 1\modu p$, hence $q$ is a prime ideal of residual degree $1$.
$q_1 \Z[\zeta_p]$ is a prime ideal of $O_K$. Let us denote $\mathbf q_1$ this ideal.
From Stickelberger relation
$g(\mathbf q_1)^p\Z[\zeta_p]=\mathbf q_1^{P(\sigma)}$ and so there exists $\varepsilon\in\Z[\zeta_p+\zeta_p^{-1}]^*$ such that
$g(\mathbf q_1)^p=q_1^{P(\sigma)}\times \varepsilon$ and so
\begin{displaymath}
(\frac{g(\mathbf q_1)}{\overline{g(\mathbf q_1)}})^p= (\frac{q_1}{\overline{q_1}})^{P(\sigma)}.
\end{displaymath}
From hypothesis $\frac{q_1}{\overline{q_1}}\equiv 1\modu\pi^{p+1}$ and so
$(\frac{g(\mathbf q_1}{\overline{g(\mathbf q_1}})^p\equiv 1\modu \pi^{p+1}$.
Therefore $v_\pi(g(\mathbf q_1)^p+1)>p$ if not we should have $v_\pi(\frac{g(\mathbf q_1}{\overline{g(\mathbf q_1}})^p+1)=p$.
From lemma \ref{l602061} p. \pageref{l602061} it follows that
$p^{(q-1)/p}\equiv 1\modu q$.
\end{proof}
\end{thm}
%
\paragraph{Acknowledgments:}
I thank Professor Preda Mihailescu for helpful   e-mail dialogues  and suggestions and for
error detections in intermediate versions of this paper.
%

%

%
Roland Qu\^eme

13 avenue du ch\^ateau d'eau

31490 Brax

France

mailto: roland.queme@wanadoo.fr

home page: http://roland.queme.free.fr/

************************************

V02 - MSC Classification : 11R18;  11R29

************************************
\end{document}